\documentclass[12pt]{article}
\usepackage{fancybox}

\everymath{\displaystyle}
\usepackage{amsmath,amsfonts,amssymb}
\usepackage{amsthm}
\theoremstyle{plain}
\usepackage{graphicx,epsfig, color}
\usepackage{mathtools}

\def\ss{\par\smallskip}
\def\ms{\par \medskip}
\def\bs{\par \bigskip}
\def\set#1{\{\,#1\,\}}

\def\inverse{^{-1}}

\def\bigvert{\ \big\vert\ }
\def\equaldef{ \buildrel\hbox{\scriptsize  def}\over= }

\def\half{{\frac{1}{2}}}

\def\qed{{\hfil \mbox{$\Box$}}}
\def\and{\ \mbox{and} \ }
\def\be{\begin{enumerate}}
\def\ee{\end{enumerate}}
\def\bi{\begin{itemize}}
\def\ei{\end{itemize}}
\def\beqn{\begin{eqnarray*}}
\def\eeqn{\end{eqnarray*}}
\def\n{\mathbb N}

\def\F{\mathbb F}
\def\r{\mathbb R}
\def\c{\mathbb C}
\def\c{\mathbb C}

\def\half{\frac{1}{2}}

\def\bc{\begin{center}}
\def\ec{\end{center}}
\def\ef{\end{frame}}
\def\ed{\end{document}}
\def\gFpair{\big(g(x),\, F(x)\big)}
\def\om{\omega}
\def\order{\text{order}}
\def\wg{\widehat{g}}
\def\ord{\text{ord}}
\def\ellcm{\ell cm}
\def\SigF{{\textstyle \Sigma_F}}
\begin{document}

\newtheorem{thm}{Theorem}
\newtheorem*{thm*}{Theorem}
\newtheorem{ex}{Example}[section]
\newtheorem{defn}[ex]{Definition}
\newtheorem{cor}[ex]{Corollary}
\newtheorem{prop}[ex]{Proposition}
\newtheorem{lem}[ex]{Lemma}
\newtheorem{noat}[ex]{Note}
\newtheorem{rmk}[ex]{Remark}
\newtheorem{notation}[ex]{Notation}
\newtheorem*{qstn*}{Open question}
\bc
{\bf\large Elements of Finite Order in the Riordan Group
\ms
and Their Eigenvectors}
\ms
{\sc  Marshall M. Cohen}
\bs
{\small {\it Department of Mathematics, Morgan State University, Baltimore, MD}\\
Email: marshall.cohen@morgan.edu}
\bs

\ec
\ms

\hrulefill
\ms
\bc
\parbox{6in}{ {\bf Abstract.}\sl \ We consider elements of finite order in the Riordan group $\cal R$ over a field $\F$ of characteristic $0$.  Viewing  $\cal R$ as a semi-direct product of groups of formal power series, we solve, for all $n \geq 2$,  two foundational questions posed by L. Shapiro [16] for the case $n = 2$ (``involutions''): 
Theorem 1 states, 
given  $F(x)$ of finite compositional order, 
exactly which $g(x)$ make $\big(g(x), F(x))$ a Riordan element of order $n$.  Theorem 2 classifies finite-order Riordan group elements up to conjugation in $\cal R$.  Viewing $\cal R$  as a group of infinite lower triangular matrices,  we interpret Theorem 1 in terms of existence of eigenvectors and Theorem 2 as a normal form for finite order Riordan arrays under similarity. These lead to Theorem 3,  a formula for all eigenvectors of finite order Riordan arrays; and we show how this can lead to interesting combinatorial identities. We then relate our work to papers [4], [5] of Cheon and Kim which motivated this paper and we solve the Open question in Section 2 of [5].   Finally, this circle of ideas gives a new proof of C. Marshall's theorem ([9], [10]), which finds the unique $F(x)$, given bi-invertible $g(x)$, such that $\big(g(x), F(x))$ is an involution.}
\ec
\hrulefill
{\small MSC: 05A15, 20Hxx\\
Key terms: Riordan group, finite order, conjugation, eigenvector, combinatorial identities.}
\section {Introduction}
\ms
\subsection*{Background}

The {\sl Riordan group} over a field $\F$ is a group ${\cal R} = {\cal R}(\F)$ of infinite lower triangle matrices,  which form a vast generalization of Pascal's Triangle. Each element $[d_{n, j}]_{n, j \geq 0}$ of $\cal R$ is determined by a pair $\gFpair$ of generating functions, 
\[g(x) = g_0 + g_1x + \ldots,\ F(x) = f_1x + f_2x^2 + \ldots,\ \ g_0\neq 0, f_1 \neq 0,\]  where
the $j^{\rm th }$ column $(j = 0, 1, 2,  \ldots)$ is the sequence with generating function
$g(x)\cdot F(x)^j. $ Thus (letting $[x^n]\,h(x)$ denote the coefficient of $x^n$ in $h(x)$)
\[d_{n, j} = [x^n]\,g(x)\cdot F(x)^j.\]
\ss
We get Pascal's triangle, $P = [d_{n, j}] = \left[{n \choose j}\right] \ (\text{with}\ d_{n, j}= 0 \ \text{if}\  j > n)$,  if we set
\[P = \gFpair = \left(\frac{1}{1 - x}, \ \frac{x}{1 - x}\right).\]
The Riordan group was introduced in 1991 by Shapiro, Getu, Woan and Woodson [15] for combinatorial purposes. They give examples of uses as a tool in enumeration (the number of paths traversed in $n$ steps under given rules), combinatorial identities and inverse relations and such combinatorial applications have been much studied.  (See for example [1], [7], [8] and [17].)

\subsection*{The Setup}

Each element of $\cal R$ is determined by a pair $\gFpair$ of generating functions and the group is isomorphic  to the {\it semi-direct product} of two groups of formal power series.
We move seamlessly between the semi-direct product and matrix viewpoints and write, for example (the most elementary but very important example),
\[   (g_0, \omega x) = \left[\begin{array}{ccccc}
       g_0&0&0&0&\cdots\\
       0&g_0\omega&0&0&\cdots\\
       0& 0 &g_0\omega^2&0&\cdots\\
       0& 0 & 0 &g_0\omega^3& \cdots\\
      \vdots &\vdots&\vdots&0&\ddots
      \end{array}\right], \ \ g_0, \omega \in \F\setminus \{0\}.
\]
\ms
  In this paper, we  resolve some of the foundational algebraic problems concerning Riordan groups (Theorems 1 and 2) and we interpret these in terms of matrices to derive (Theorem 3) the complete determination of the  eigenvectors of a finite-order Riordan array.
\ms

\bs
We fix a field $\F$  of characteristic zero and  let $F[[x]]$  denote the set of all  {\bf formal power series} $g(x)$, 
\[g(x) = g_0 + g_1x + g_2x^2 \cdots + g_nx^n + \cdots \,, \  \   (g_n \in \F).\]
\ms
We consider the following three groups. (Technical background will be given in Section 2.)
\be
\item $\F_0[[x]] = \set{g(x) \in \F[[x]] \bigvert g_0 \neq 0}$ is a group under {\sl multiplication}.\\
\ms
\item $\F_1[[x]] = \left\{G(x) \in \F[[x]] \bigvert G(x) = g_1x + g_2x^2 + \cdots, \text{with} \ g_1\neq 0\right\}$ is a group under \\ {\sl composition}. \quad  We denote $\big(\underbrace{G\circ G\circ \cdots \circ G}_{n\ \text{terms}}\big)(x) = G^{(n)}(x)$\\
\ms
\item $\F_0[[x]]\times \F_1[[x]]$ becomes the Riordan group ${\cal R}(\F)$  if we define the operation as
\[\big(g(x), F(x)\big)\,\big(h(x), K(x)\big) = \big(g(x)\cdot h\big(F(x)\big), K(F(x))\big)\]
\ms
\item It is also useful to have specific notation for the subset $\F_+[[x]] $ of $\F[[x]]$ given by:
\ms
$\F_+[[x]] = \left\{G(x) \in \F[[x]] \bigvert G(x) = g_kx^k + g_{k + 1}x^{k + 1} + \cdots, \text{with} \ k \geq 1 \and g_k \neq 0\right\}$
\ms
{\bf Notation}:
\ms
We denote series in $\F_0[[x]]$ with lower case letters $g(x), h(x)$, etc.
\ms
We denote series in $\F_+[[x]]$  with capital letters $F(x), G(x)$, etc.
\ms
A non-constant series $g(x)\in \F_0[[x]]$  may be written uniquely as $g_0 + G(x)$ with $0\neq g_0\in \F$.
\ee
\begin{defn} \mbox{}
\ms
The {\bf order} of an element $\gamma \in \Gamma$ (in this paper $\Gamma$ will be one of the groups $\F\setminus \{0\},\ \F_0[[x]], \F_1[[x]], \ {\cal R}(\F)$) is the least possible integer $n$ such that $\gamma^n = \text{\rm (the identity element of $\Gamma$).}$ If no such integer exists, we say that $\gamma$ has {\bf infinite order}. If $\gamma$ has order two, $\gamma$ is called an {\bf involution} in $\Gamma$.
\end{defn}

\begin{rmk} \mbox{}
\ms
If $\F$ is a subfield of $\r$ then it will be immediate from the definitions of the operations in our groups (Lemma 2.1)  that the only non-trivial elements of finite order in these groups are involutions, since this  is true in the multiplicative group $\r\setminus \{0\}$.  This fact is relevant when dealing with generating functions of integer (e.g., counting) sequences as is central in Combinatorics.  Theorem 1 below gives a significant new result even when restricted to involutions.  On the other hand, power series over the complex field $\F = \c$ are  central in mathematics and our main theorems answer basic questions in this realm for all $n\geq 2$. \end{rmk}
\subsection*{Main Theorems} 
\begin{thm} \mbox{}
\ss
Suppose that $g(x) = g_0 + g_1x + g_2x^2 + \cdots \, \in \F_0[[x]]$, where $g_0$ has finite order $(=\ord(g_0))$ in $\F\setminus \{0\}$ and $F(x) = \om x + F_2x^2 + \cdots $ is an element of finite (compositional) order in $\F_1[[x]]$.
\ Then
\ms
 $\ \gFpair$ has finite order $n$ in the Riordan group \ \ 
$\iff$ 

\be
\item $n = \ell cm\big(\text{ord}(g_0), \text{ord}(F(x))\big) = \ellcm\big(\text{ord}(g_0), \text{ord}(\om)\big) $
\item [{\bf and}]
\item
there exists an element $h(x) \in \F_0[[x]]$ such that $g(x) = g_0\cdot\frac{h(x)}{h\big(F(x)\big)}$.
\ee
\bs\ss
Indeed, if \ $\wg(x) \equaldef \frac{1}{g_0}\cdot g(x)$ \ and \ \ $b \equaldef \ord\big(F(x)\big)$
  then   $g(x) = g_0\cdot\frac{h(x)}{h\big(F(x)\big)}$ if we set
\[h(x) = \Big(\wg^\frac{1}{b}(x)\Big)^{b - 1}\Big(\wg^\frac{1}{b}\big(F(x)\big)\Big)^{b - 2}\Big(\wg^\frac{1}{b}\big(F^{(2)}(x)\big)\Big)^{b - 3} \cdots \ \Big(\wg^\frac{1}{b}\big(F^{(b - 2)}(x)\big)\Big)\]
\end{thm}
\begin{rmk}\ \ 
  If $m\in \n$ then $\wg(x) = \frac{1}{g_0}\cdot g(x) =  (1 + b_1x +\cdots\,)$ has a unique $m^{\rm th}$ root of the form $(1 + c_1x + \cdots\,)$ by  [11], Theorem 3. It is this $m^{\rm th}$ root which is denoted by $\wg^{\frac{1}{m}}(x)$. It is proved in [11] that all the usual laws of roots and exponents apply for these formal power series. In our proof of Theorem 1 we could replace $\wg^{\frac{1}{b}}(x)$ by $\alpha\cdot \wg^{\frac{1}{b}}(x)$ in the formula for $h(x)$, where $\alpha$ is any $b^{\rm th}$ root of unity.
\end{rmk}
\begin{rmk} If $\gFpair$ is an {\bf involution} (i.e., n = 2) then the proof gives $h(x) = 
\Big(\frac{1}{g_0}\cdot g(x)\Big)^\half$.\\
\end{rmk}

\begin{thm}[The Conjugacy Theorem] \mbox{}
\be
\item [A.] Suppose that $g(x) = g_0 + g_1x + g_2x^2 + \cdots$ and $F(x) = \omega x + F_2x^2 + \cdots$.
\ms
If $\gFpair$ has finite order in the Riordan group then $\gFpair$ is conjugate in ${\cal R}$ to $\big(g_0, \omega x\big)$. \\
\bs
  Indeed if  $h(x)\in \F_0[[x]]$ with $g(x) = g_0\cdot\frac{h(x)}{h\big(F(x)\big)}$, as in Theorem 1,  if $\ord(F(x) = b$ and ${\textstyle \Sigma_F = \frac{1}{b}\sum_{j=1}^b \omega^{j}F^{(b - j)}(x)}$ then
\[\Big(h(x), \ \Sigma_F\Big)\inverse\,\Big(g(x), F(x)\Big) \Big(h(x), \Sigma_F 
 \Big) = \big(g_0, \omega x\big)\]
 \item [B.] Suppose that $\big(g(x), F(x)\big) \and  \big(h(x), K(x)\big)$ are Riordan elements of finite order with \[{\textstyle g(x) = \sum_i g_ix^i, \quad h(x) = \sum_i h_ix^i, \quad F(x) = \sum_j f_jx^j 
 \ \  \and \  K(x) = \sum_j k_j x^j.}\]
 \ms
 Then \ \ $\Big(g(x), F(x)\Big) \sim  \Big(h(x), K(x)\Big)$ \ \ ({\it i.e.,} they are conjugate in $\cal R$)
 \ms
 \mbox{\hspace{.4in}} $\iff\ \  \ g_0 = h_0\  \and \ f_1 = k_1$.
 \ee
\end{thm}
\ms

\begin{rmk}
 If $\gFpair$ is an {\bf involution} and $\gFpair \neq (-1, x)$ then the Conjugacy Theorem, Part A, states that 
\[\Big(\big(\frac{1}{g_0}\cdot g(x)\big)^\half, \ \frac{1}{2}\big(x - F(x)\big)\Big)\inverse\,\Big(g(x), F(x)\Big)\, \Big(\big(\frac{1}{g_0}\cdot g(x)\big)^\half, \ \frac{1}{2}\big(x - F(x)\big)\Big) = \big(1, -x\big).\]
\end{rmk}
\bs
Our third main theorem,
\ms
{\bf Theorem 3. (The general form of an eigenvector),}
\ms
will be stated and proved in Section 5.   It gives a complete formula for the eigenvectors of a finite order element $\gFpair$.
\begin{rmk}  \mbox{}
\ms
While our Main Theorems are logically independent of their work,
we were much influenced by the papers ([4], [5]) of Cheon and Kim.  We discuss connections to their papers in Section 6 and, in particular, we solve the Open question in Section 2 of [5]. In Section 7 we use an idea from [4] to give a new (but less elegant) proof of C. Marshall's Theorem ([9], [10]).

\end{rmk}

{\bf Acknowledgments.} I wish to acknowledge my colleagues at Morgan State University, Drs. Leon Woodson, Asamoah Nkwanta, Xiao-Xiang Gan and Rodney Kerby, and Dr. Lou Shapiro of Howard University, who introduced me to these matters and provided the collegial and stimulating atmosphere in which to work on them. And I am grateful to Dr. Gi-Sang Cheon for pointing out that Theorem 1 could be restated in terms of eigenvectors.
\section {Basics of Formal Power Series and the Riordan Group}
We summarize the initial facts  about formal power series  and the Riordan group with Lemmas 2.1 - 2.3. See [6] for proofs, references and historical background on formal power series, [1] and [15] for introduction to the Riordan group. In 2.4 - 2.7 we develop some basic facts about elements of finite order in $\cal R$ which we will need in the next section.
\ms
Our notational convention is to follow the pattern
\[g(x) = \sum_{n = 0}^\infty g_nx^n\ \in \F_0[[x]] ,\ \ \ F(x) = \sum_{n=1}^\infty f_nx^n \in \F_+[[x]],\ \ {\rm etc.,}\]
\begin{lem} \mbox{}
\bi 
\item $\F_0[[x]]$ is an abelian group under multiplication:
\bi
\item $g(x)\cdot h(x) \, \equaldef \, g_0h_0 + (g_0h_1 + g_1h_0)x \, +\,  \cdots \, + \, \left({\textstyle \sum\limits_{k=0}^n}(g_kh_{n - k}\right)x^n + \cdots\, $ 
\item The identity element is $g(x) = 1$.
\item $g(x)\inverse \equaldef \frac{1}{g(x)} = \frac{1}{g_0} - \frac{g_1}{g_0^2}\,x + \cdots$\ \ (found by solving inductively for coeffients).\\
Or, writing $g(x) = g_0\big(1 + G(x)\big)$,  note: \  $\frac{1}{g(x)} =\frac{1}{g_0}\Big(1 - G(x) + G(x)^2 - \cdots )$ \ (geometric series).
\ei
\item $\F_1[[x]]$ is a group under composition (= {\rm ``Substitution}'')
\bi
\item $(F\circ G)(x) = F\big(G(x)\big) \equaldef f_1\cdot G(x) + f_2\cdot G(x)^2 + f_3G(x)^3 + \ldots$ 
\ss
\mbox{\hspace{.6in}} $ = f_1g_1x + (f_1g_2 + f_2g_1^2)x^2 + \cdots $
\ms
\item The identity is $\text{\rm id}(x) =x$.
\ms
\item The $n^{\rm th}$ power of $F$ is denoted by $F^{(n)}(x) = F\big(F(\cdots F(x)\cdots\big))$.
 \ss
  \mbox{\hspace{2.1in}}$F^{(n)}(x) = f_1^nx + \cdots $
\ms
\item The inverse $F^{(-1)}(x)$of $F(x)$ in $\F_1[[x]]$ is denoted by $\overline{F}(x) = F^{(-1)}(x)$:
\ss
\centerline{$\overline{F}(x) = \frac{1}{f_1}x - \frac{f_2}{f_1^3}x^2 + \cdots $\ \ \  (solve $F(G(x)) = x$ inductively for the coefficients $g_j$).}
\ms
(For a general formula, see the {\sl Lagrange Inversion Formula}, [19].)
\ei
\item The set $\F_0[[x]] \times \F_1[[x]]$ becomes a group -- the Riordan group 
${\cal R} = {\cal R}(\F)$ -- if we define
\[\big(g(x), F(x)\big)\big(h(x), K(x)\big) = \big(g(x)\cdot h(F(x)), \, K(F(x)).
\qquad\]
\bi
\item The identity element in ${\cal R}$ is $\big(1, x)$ 
\ms
\item Inverses are given by \ \ $\Big(g(x), \ F(x)\Big)\inverse = 
\left(\frac{1}{g(\overline{F}(x))},\ \overline{F}(x)\right).$
\ms
\item If $n\in \n$ then $\big(g(x), F(x)\big)^n = \Big(g(x)\cdot g(F(x)) \cdots g(F^{(n - 1)})(x), \ F^{(n)}(x)\Big) \\
\mbox{} \hspace{1.9 in} \ = \big(g_0^n + \cdots\ ,\  f_1^n\,x + \cdots \,\big).$
\ei
\ei
\end{lem} 
\begin{lem}\mbox{}\\
If \quad $F(x) = \omega x + f_2x^2 + f_3x^3 + \ldots$ has finite (compositional) order $n$ in $\F_1[[x]]$,
\ms
 then $\omega$ has order $n$ in $\F\setminus \{0\}$.\ \  (i.e., $\omega$ is a primitive $n^{\rm th}$ root of unity.) $\Box$
\end{lem}
\bs
The author first learned the following Conjugation Lemma from a much more general theorem (Theorem 8, page 19) in [2] . It appears explicitly in [5] and [14] and appears with different conjugator than $\SigF$ in [3].
\begin{lem} {\rm (Conjugation of finite-order formal power series: $F(x) \sim \om x$)}
\ms
Suppose that $F(x) = \om x + f_2x^2 + f_3x^3 + \cdots $ has finite order $n$ in $\F_1[[x]]$. 
Let \[{\textstyle\Sigma_F} = \frac{1}{n}\sum_{j = 1}^n \om^jF^{(n - j)}(x).\] 
\ms
Then $\SigF$ conjugates $F(x)$ to $\om x$ in $\F_1[[x]]$: \ $\big({\textstyle\SigF}\circ F\circ {\textstyle\SigF}\inverse\big)(x) = \ell_\om(x). \ \  \Box$ 
\end{lem}
\bs
The following simple examples of  Riordan elements of order $6$ in ${\cal R}(\c) $ illustrate the next Proposition.
\begin{ex} 
$(a, bx)(c, dx) = (ac, (db)x),\ a, b, c, d \in\F$. \\  
Thus all of the following elements of ${\cal R}\big(\c\big)$
have order 6:
\[\big(1, e^{\frac{2\pi i}{6}}z\big),\quad   \big(e^{\frac{2\pi i}{3}}, -z\big), \quad \big(e^{\frac{2\pi i}{6}}, -z\big), \quad \big(e^{\frac{2\pi i}{6}}, e^{\frac{2\pi i}{6}}z\big) \]
\end{ex}

The following Proposition is central to the study of Riordan elements of finite order $n \geq 2$ and to the proofs of our Main Theorems.
\begin{prop} \mbox{}\\
Let $g(x) = g_0 + g_1x + g_2x^2 + \cdots \, \in \F_0[[x]]$ and $F(x)\in \F_1[[x]]$. Suppose that 
\bi
\item $g_0$ has order $a$ in the multiplicative group $\F\setminus \{0\}$.
\item $F(x)$ has (compositional) order $b$ in the  group $\F_1[[x]]$
\item $\gFpair$ has finite order $n$ in the Riordan group.
\ei
\ms
Then $n = \ell cm\big(a,\,  b\big)$.
\end{prop}
\ms
{\bf Proof:}
\ms
A straightforward inductive argument, using the definition of multiplication in ${\cal R}$ yields the fact that $\gFpair$ has finite order $n$ iff $n$ is the least positive integer such that
\[\bullet\ \ g(x)\cdot g(F(x)) \cdots \cdot g\big(F^{(n - 1}(x)\big) = 1 \quad \  \text{\bf and}\quad \ \  \bullet \ \ F^{(n)}(x) = x. \qquad \]
\ss
Since $g(x)\cdot g(F(x)) \cdots \cdot g\big(F^{(n - 1)}(x)\big) = (g_0^n + \cdots)$ we see that $g_0^n = 1$ so that $a = \order(g_0)$ divides $n$. Also $F^{(n)}(x) = x$ implies  that $b$ divides $n$.  Hence $\ell cm(a, b)$ divides $n$.
\ms
Let $n = kb$. Then
\beqn
(\big(1, x\big) = \gFpair^n & = &  \Big(g(x)\cdot g(F(x) \cdots g(F^{(kb-1)}(x), \ F^{(kb)}(x)\Big) \qquad \qquad \\
&&\\
& = & \Big(\big(g(x)\cdot g(F(x) \cdots g(F^{(b-1)}(x))\big)^k, \ x\Big)\ \text{(because}   F^{(jb + r)}(x) = F^{(r)}(x)).\\
\eeqn

This implies that \ $g(x)\cdot g(F(x))\cdots g(F^{(b-1)}(x)) = (g_0^b + \cdots)$ is a constant (a $k^{\rm th}$ root of unity in $\F$). Hence
\[g(x)\cdot g(F(x) \cdots g(F^{(b-1)}(x) = g_0^b.\]
Now let us set $\ell cm = \ell cm(a, b) = qb$.  Thus 
\ $1 = g_0^a = g_0^{\ell cm} = \big(g_0^b\big)^q$\ and \ $x = F^{(b)}(x) = F^{\ell cm}(x)$.
Therefore 
\beqn
\gFpair^{\ell cm} & = & \Big(\big(g(x)\cdot g(F(x) \cdots g(F^{(b-1)}(x)\big)^q, \ F^{\ell cm}(x)\Big)\\
&&\\
& = & \Big(\big(g_0^b\big)^q,\ x)\\
&&\\
& = & \big( 1, x\big)
\eeqn
Hence $n = \order\gFpair$ divides $\ell cm(a, b)$  and, as seen at the outset, $\ell cm(a, b)$ divides $n
$. Thus $n = \ell cm(a, b)$, as claimed. \quad $\Box$
\ms
\begin{cor} \mbox{}
\ss
Suppose that $g(x) = 1 + g_2x^2 + \cdots$ and that $\gFpair$ has finite order.
\ss
Then $\order\gFpair = \order(F(x))$ \quad $\Box$
\end{cor}
\ss
\begin{cor} \mbox{}
\ss
Suppose that  $g(x) = g_0\cdot k(x)$ where $g_0\in \F \setminus \{0\}$ has finite order and $k(x)\in \F_0[[x]]$.  Let $F(x) \in \F_1[[x]]$.
\be
\item $\gFpair$ has finite order iff $\big(k(x),\, F(x)\big)$ has finite order.
\item When these have finite order, 
\item []$\bullet \ \ord\gFpair  = \ell cm\big(\ord(g_0k_0), \ord(F(x))\big) \qquad \bullet \ 
\ord\big(k(x), F(x)\big) = \ell cm\big(\ord(k_0), \ord(F(x)\big)$
\ee
\end{cor}
{\bf Proof:} 
\be
\item If $\gFpair$ has finite order $n$ and if $\ell = \ell cm(n, \order(g_0)$ then 
\beqn
(1, x) & = & \gFpair^n = \gFpair^\ell\\
& = & \Big(g(x)\cdot g(F(x)) \cdots g(F^{(\ell -1)}(x)), \ F^{(\ell)}(x)\Big) \\
& = & \Big(g_0^\ell \cdot k(x)\cdot k(F(x)) \cdots k(F^{(\ell - 1)}(x)), \ F^{(\ell)}(x)\Big)\\
& = & \Big(k(x)\cdot k(F(x)) \cdots k(F^{(\ell-1)}(x)), \ F^{(\ell)}(x)\Big) \\
\eeqn
Thus $\big(k(x), F(x)\big)$ has finite order dividing $\ell$.
\ms
Similarly, if $\big(k(x), F(x)\big)$ has order $m$, then $\gFpair^m =  \Big(g_0^m \cdot k(x)\cdot k(F(x)) \cdots (F^{(m - 1)})(x), F^{(m)}\Big)$.  Thus $\gFpair^{\ell cm\big(\ord(g_0), m\big)} = (1, x)$ and $\gFpair$ has finite order.

\item The evaluation of the orders of $\gFpair \and \big(k(x), F(x)\big)$ follow from Proposition 2.5. \quad $\Box$
\ee

\section{Proof of Theorem 1: $g(x) = g_0\cdot\frac{h(x)}{h(F(x))}$}
$(\Longleftarrow)$\  {\bf Proof of Sufficiency in Theorem 1}
\ms
We are given that $g(x) = g_0 + G(x)\in \F_0[[x]]$ and $F(x)\in \F_1[[x]]$ with $g_0$ of finite order and  $F(x)$ of finite order. We let  $n = \ellcm(\ord(g_0),\, \ord(F(x))$.  Also we are given $h(x)\in \F_0[[x]]$ with $g(x) = g_0\cdot\frac{h(x)}{h\big(F(x)\big)}$.
Then we have
\beqn
\gFpair^n & = &  \Big(g(x)\cdot g(F(x)) \cdots \, g(F^{(n - 1)}(x)), \ F^{(n)}(x)\Big)\\
&&\\
& = & \Big(g_0^{n}\cdot\frac{h(x)}{h\big(F(x)\big)}\frac{h(F(x))}{h\big(F^{(2)}(x)\big)}\ \cdots \frac{h(F^{(n - 1)} (x))}{h(F^{(n)}(x))},\ x\,\Big)\\
&&\\
& = & \Big(1\cdot\frac{h(x)}{h\big(F(x)\big)}\frac{h(F(x))}{h\big(F^{(2)}(x)\big)}\ \cdots \frac{h(F^{(n - 1)} (x))}{h(x)},\  x\,\Big) \ = \ \big(1, x\big)\\
\eeqn%
Therefore $\gFpair$ has finite order. We apply Corollary 2.7, setting $k(x) = \frac{h(x)}{h\big(F(x)}\big) = 1 + K(x)$,  and we see that $\ord\gFpair = \ellcm\big(\ord(g_0),\ \ord(F(x))\big)  = n$. \quad $\Box$
\bs\ms
$(\Longrightarrow)$\  {\bf Proof of Necessity in Theorem 1}  
\ms
We are given that $\gFpair$  has order $n$. The fact that $n = \ell cm\big(\text{ord}(g_0), \text{ord}(F(x))\big)$ is given by 
Proposition 2.5.  It remains to prove that $\exists \ h(x)\in \F_0[[x]]$ such that \break $g(x) = g_0\cdot\frac{h(x)}{h\big(F(x)\big)}$. 
\bs
 If $g(x)  \in \F_0[[x]]$ we set\ \ 
$\widehat{g}(x) = \frac{1}{g_0}\cdot g(x) = 1 + \frac{g_1}{g_0}\,x + \frac{g_2}{g_0}\,x^2 + \cdots $.
\ms
 From Corollary 2.7, we know that $\big(\wg(x),\, F(x)\big)$ is an element in the Riordan group of order equal to $\ellcm\big(1, \ord(F(x)\big) = \ord(F(x))$.  If we set $b = \ord(F(x))$ then we have
\[1  =  \wg(x)\wg\big(F(x)\big) \cdots \wg\big(F^{(b - 1)}(x)\big)\qquad \and \qquad  
\frac{1}{\wg(x)}  =  \wg\big(F(x)\big)\,\wg\big(F^{(2)}(x)\big)\cdots \wg\big(F^{(b - 1)}(x)\big)\quad (**)\]
Using Remark 1.3, the  series $\wg^{\frac{1}{b}}(x)$ is well-defined and for notational convenience we denote $k(x) = \wg^{\frac{1}{b}}(x)$ .  Let us define
\beqn
h(x) & \equaldef & \Big(\wg^\frac{1}{b}(x)\Big)^{b - 1}\Big(\wg^\frac{1}{b}\big(F(x)\big)\Big)^{b - 2}\Big(\wg^\frac{1}{b}\big(F^{(2)}(x)\big)\Big)^{b - 3} \cdots \ \Big(\wg^\frac{1}{b}\big(F^{(b - 2)}(x)\big)\Big)^{1}\\
&&\\
& = &  \big(k(x)\big)^{b - 1}\cdot \big(k(F(x))\big)^{b - 2}\cdot \big(k(F^{(2)}(x))\big)^{b - 3}
\cdots  \cdot\big(k(F^{(b - 2)}(x))\big).
\eeqn
Then
\beqn
g(x) & = & g_0\cdot\wg(x) \ = \ g_0\left[\,\frac{ \widehat{g}(x)^{b - 1}}{\Big(1\slash \widehat{g}(x)\Big)}\,\right]^{\frac{1}{b}}\\
&&\\
&&\\
& = & g_0\left[\,\frac{ \widehat{g}(x)^{b - 1}}{(\widehat{g}\big(F(x)\big)\cdot\widehat{g}\big(F^{(2)}(x)\big)\cdot \cdots \cdot \widehat{g}\big(F^{(b-1)}(x)\big)}\,\right]^{\frac{1}{b}} \quad \text{by equation $(**)$ above}\\
\\
&&\\
& = &  g_0\cdot \frac{ k(x)^{b-1}}{\big(k(F(x)\big)\big(k(F^{(2)}(x))\big)\big(k(F^{(3)}(x))\big) \cdots \cdot\big(k(F^{(b - 1)}(x))\big)}\\
&&\\
 & = & g_0\cdot\frac{\big(k(x)\big)^{b - 1}}{\big(k(F(x))\big)^{b - 1}}
\cdot\frac{\big(k(F(x))\big)^{b - 2}}{\big(k(F^{(2)}(x))\big)^{b - 2}}
\cdot \frac{\big(k(F^{(2)}(x))\big)^{b - 3}}{\big(k(F^{(3)}(x))\big)^{b - 3}}
\cdots  \cdot\frac{\big(k(F^{(b - 2)}(x))\big)}{\big(k(F^{(b - 1)}(x))\big)}\\
&&\\
& = &  g_0\cdot\frac{h(x)}{h(F(x)}\\
\eeqn
{This completes the proof of  Theorem 1. }\quad $\Box$
\section{Proof of Theorem 2 (The Conjugacy Theorem)}
\be
\item [A.] From Theorem 1, there exists $h(x)\in \F_0[[x]]$ such that $g(x) = g_0\cdot\frac{h(x)}{h(F(x)}$. \\ Let $\Sigma_F = \sum_{j=1}^b \omega^jF^{(b - j)}(x) \in \F_1[[x]]$, as in Lemma 2.3.  Denote $\ell_\om(x) = \om x$. \ Then
\beqn
\gFpair\,\big(h(x), \Sigma_F)\, & = & \Big(g(x)\cdot h(F(x)), \, (\Sigma_F\circ F)(x)\Big)\\
&&\\
& = & \Big(g_0\cdot\frac{h(x)}{h(F(x)}\cdot h(F(x)), \, (\Sigma_F\circ F)(x)\Big)\\
&&\\
& = &\Big(g_0\cdot h(x), \, (\ell_\om\circ \Sigma_F)(x)\Big), \quad \text{by Lemma 2.3,}\\
&& \\
& = & \big(h(x), \Sigma_F(x)\big)\, \big(g_0, \ell_\om(x)\big)\\
&=& \big(h(x), \Sigma_F(x)\big)\, \big(g_0, \om x\big)\\
\eeqn

Therefore \qquad $\big(h(x), \Sigma_F\big)\inverse\,\gFpair\,\big(h(x), \Sigma_F\big) = 
\big(g_0, \om x)$, \quad as claimed. \ \ \ $\Box$
\ms
\item [B.] If $g_0 = h_0 \and f_1 = k_1$ then, from part A., 
\[\big(g(x), F(x)\big) \sim \big(g_0, f_1\big) = \big(h_0, k_1\big) \sim \big(h(x), K(x)\big).\]
Therefore \ $\big(g(x), F(x)\big) \sim \big(h(x), K(x)\big)$.
\ms
Conversely suppose \  $\big(g(x), F(x)\big) \sim \big(h(x), K(x)\big)$.  We observe that if we conjugate $\big(g(x), F(x)\big)$ by an element $\big(a(x), B(x)\big)$, then the first coefficients of the two terms of the pair both remain unchanged:
\beqn
 &&\Big(a(x), B(x)\Big)\Big(g(x), F(x)\Big)\Big(a(x), B(x)\Big)\inverse \mbox{\hspace{2.5in}}\\
 & = & \Big(a(x)\cdot g(B(x)),\ F(B(x))\Big)\Big(\frac{1}{(a\big(\overline{B}(x)\big)},\  \overline{B}(x)\Big)\\
& = & \Big(a(x)\cdot g(B(x)\cdot\frac{1}{(a(\overline{B}\Big(F(B(x))\Big)},\ \overline{B}(F(B))\Big)\\
& = & \Big((a_0 +\cdots)(g_0 +\cdots)\frac{1}{(a_0 + \cdots)}\,,\ (\frac{1}{b_1}\cdot f_1\cdot b_1)x + \cdots \Big)\\
& = &  \Big((g_0 +\cdots,\ f_1\cdot x + \cdots \Big)
\eeqn
Therefore $g_0 = h_0 \and f_1 = k_1$.\quad $\Box$
\ee
\bs
\section {Eigenvectors of Riordan Elements of Finite Order}
\ms
Our algebraic theorem giving classification of finite order Riordan elements up to conjugation becomes, from the point of view of linear algebra, a normal form for finite order proper arrays under similarity.  We use the fact that similar matrices have corresponding eigenvectors: if $\vec{h}$ is an eigenvector of the matrix $B$  and if $A = PBP\inverse$ then $\vec{v} = P(\vec{h})$ is an eigenvector of $A$ with the same eigenvalue.
\ms
We are grateful to Dr. Gi-Sang Cheon for the following observation (see also [1], Section 6.6), which connects our main theorem to eigenvectors of Riordan arrays.
\ms
\begin{lem}
If $\gFpair \in {\cal R}$ and $\vec{h} = (h_0, h_1, \ldots\ )^T$ \ with generating function $h(x)$ and $h_0\neq 0$\
then $\vec{h}$ is an eigenvector of $\gFpair$ with eigenvector $g_0$ if and only if
\beqn
\gFpair\cdot h(x) & = & g(x)\cdot h\big(F(x)\big) = g_0\cdot h(x).\\
\iff g(x) & = &  g_0\cdot \frac{h(x)}{h(F(x)}\\
\eeqn
\end{lem}
{\bf Note:}\ \
If $h_0 = 0$ we would only write this in the form\ $g(x)\cdot h(F(x)) = g_0\cdot h(x)$ since $h(F(x))$ would not have a multiplicative inverse.
\bs
{\bf Proof:}
The opening discussion of the Riordan group [15], concerning how multiplication of Riordan arrays is defined, is based on the fact (now known as {\it The Fundamental Theorem of Riordan Groups}) that the multiplication of a Riordan matrix $\gFpair$ by a vector $\vec{h}$ gives
\[\gFpair\cdot\vec{h} \ \ \ \text{has generating function}\ \  \  g(x)\cdot h\big(F(x)\big). \qquad \qed\]
\bs

Motivated by Theorem 1, we fix $h(x)$ in our current discussion:
\begin{notation}
\[h(x) \equaldef \Big(\wg^\frac{1}{b}(x)\Big)^{b - 1}\Big(\wg^\frac{1}{b}\big(F(x)\big)\Big)^{b - 2}\Big(\wg^\frac{1}{b}\big(F^{(2)}(x)\big)\Big)^{b - 3} \cdots \ \Big(\wg^\frac{1}{b}\big(F^{(b - 2)}(x)\big)\Big)\]
where \ $\wg(x) \equaldef \frac{1}{g_0}\cdot g(x)$ \ and \ $b = \text{ord}\big(F(x)\big)$
\end{notation}
\begin{cor}\ {\bf (Corollary to Theorem 1)}
\ms
\be
\item [(a)] If $\gFpair \in {\cal R}$ has finite order then $h(x)$ is the generating function of an eigenvector of $\gFpair$ with eigenvalue $g_0$.
\item [(b)] In particular, if $\gFpair$ is an involution in $\cal R$ then
\[h(x) = \sqrt{\,\wg(x)}\]
is the generating function of an eigenvector of $\gFpair$ with eigenvalue $g_0$.\quad $\Box$ 
\ee
 \end{cor}
 \bs
 In fact, we can find all eigenvectors of any Riordan element $\gFpair$ of finite order. Consider Theorem 2 (The Conjugacy Theorem) as a normal form for similarity of matrices, Riordan arrays of finite order.  We have
 \[\left[\begin{array}{ccccc}
       g_0&0&0&0&\cdots\\
       g_1&g_0\omega&0&0&\cdots\\
       g_2& (g_0f_2 + g_1\omega) &g_0\omega^2&0&\cdots\\
       g_3& \cdot & \cdot &g_0\omega^3&\cdots\\
      \vdots &\vdots&\vdots&&\ddots
      \end{array}\right] \ = \ 
\gFpair \ {\large \sim}\   (g_0, \omega x) = \left[\begin{array}{ccccc}
       g_0&0&0&0&\cdots\\
       0 & g_0\omega&0&0&\cdots\\
       0 & 0 &g_0\omega^2&0&\cdots\\
       0 & 0 & 0 &g_0\omega^3&\cdots\\
      \vdots &\vdots&\vdots&&\ddots
      \end{array}\right]\]
 \ms
\begin{thm} {\bf (The general form of an eigenvector)}
\ms
Suppose that $\gFpair \in {\cal R}$ has finite order and that $\text{the compositional order of $F(x)$)}\ = \ord(\omega) = b > 1$. Set \ $\wg(x) \equaldef \frac{1}{g_0}\cdot g(x)$. and  let
\beqn
 h(x) & \equaldef & \Big(\wg^\frac{1}{b}(x)\Big)^{b - 1}\Big(\wg^\frac{1}{b}\big(F(x)\big)\Big)^{b - 2}\Big(\wg^\frac{1}{b}\big(F^{(2)}(x)\big)\Big)^{b - 3} \cdots \ \Big(\wg^\frac{1}{b}\big(F^{(b - 2)}(x)\big)\Big)\\
 \Sigma_F & = &\frac{1}{b}\sum_{j=1}^b \omega^{j}F^{(b - j)}(x).
\eeqn
Then 
\ms
\bi
\item  $v(x)$ is the generating function of an eigenvector of $\gFpair$ if and only if
\[v(x) = h(x)\cdot\big(\theta\circ\Sigma_F\big)(x)\ = \big(h_0\theta_k\,x^k + \cdots\big)\]
$\text{where}\ 
\text {for some integer}\  k\geq 0,\ \theta(x) = \sum_{j = 0}^\infty \theta_{k + jb}x^{k + jb} \ \text{with}
\ \ \theta_k \neq 0.$\ 
\ss
 \item
the eigenvalue for $v(x)$ is $g_{0}\omega^k$. 
\ei
\end{thm}
\ms
\begin{cor}
Suppose that $\gFpair$ is an involution with $F(x) = -x + f_2x^2 + \cdots$. 
\ms
Then 
$v(x)$ is the generating function of an eigenvector of $\gFpair$ if and only if 
\[  v(x) = \sqrt{\wg(x)}\cdot\theta\big(x - F(x)\big). \]
where  $\theta(x) = \sum_{j = 0}^\infty \theta_{k+2j}x^{k +2j}, \ \ \theta_k\neq 0$,  is a series which is either even or odd,
\end{cor}
\ms
{\bf Proof of Theorem 3:} 
\ms
As above, we let  $ h(x) \equaldef \Big(\wg^\frac{1}{b}(x)\Big)^{b - 1}\Big(\wg^\frac{1}{b}\big(F(x)\big)\Big)^{b - 2}\Big(\wg^\frac{1}{b}\big(F^{(2)}(x)\big)\Big)^{b - 3} \cdots \ \Big(\wg^\frac{1}{b}\big(F^{(b - 2)}(x)\big)\Big).$
\be
\item From Theorem 2,
\[\Big(g(x), F(x)\Big) = \Big(h(x), \ \Sigma_F\Big)\,\big(g_0, \omega x\big)\,\Big(h(x), \Sigma_F \Big)\inverse. \]
 \item $\theta(x) = \sum_{m = k}^\infty \theta_mx^m\ (k\geq0,\theta_k\neq 0)$ is the generating function of an eigenvector \break \ss   of $\big(g_0, \omega x\big)$
 \beqn
 & \iff & \big(g_0, \omega x\big)\cdot \theta(x) = \lambda \theta(x) \hspace{3in}\\
 & \iff & g_0\cdot \theta\big(\omega x\big) = \lambda \theta(x)\\
 & \iff & \sum_{m = k}^\infty g_0\omega^m\theta_mx^m = \sum_{m = k}^\infty \lambda \theta_mx^m\\
 & \iff & (g_0\omega^m = \lambda\ \ \text{or}\ \ \theta_m = 0) \ \text {for all $m \geq k$}.
 \eeqn
 In particular, we have {\setlength\fboxsep{.1cm}\fbox{$g_0\omega^k = \lambda$}}\  so  that $\theta(x)$ gives an eigenvector of $\big(g_0, \omega x\big)$
 \beqn
 & \iff & (\omega^{m - k} = 1 \ \ \text{or}\ \ \theta_m = 0) \ \text {for all $m \geq k$} \hspace{1.6in} \\
  & \iff & \big(m = k + jb \ \  (j = 0, 1,\ldots; b = \text{ord}(\omega)) \ \ \text{or}\ \ \theta_m = 0\big) \ \text {for all $m \geq k$}\\
  & \iff & \theta(x) = \sum_{j = k}^\infty \theta_{k + jb}x^{k + jb}.
 \eeqn
 \item Since similar matrices have corresponding eigenvectors, as pointed out at the start of this section, we have, from steps 1. and 2. the fact that\\
 $v(x) = \sum_{m = 0}^\infty v_mx^m$ is an eigenvector of $\gFpair$
 \beqn
& \iff  & v(x) = \big(h(x), \Sigma_F(x)\big)\cdot\theta(x)\ \ \ \text{where}\ \theta(x) = \sum_{j = 0}^\infty \theta_{k + jb}x^{k + jb},
\ \theta_k \neq 0\\
& \iff & v(x) = h(x)\cdot\theta\big(\Sigma_F(x)\big)\ \ \ \text{where}\ \theta(x) = \sum_{j = k}^\infty \theta_{k + jb}x^{k + jb},
\ \theta_k \neq 0 \quad \Box\\
 \eeqn
\ee
\subsection* {Combinatorial Identities from Eigenvectors}
\ms
Theorem 3 leads to an interesting class of combinatorial identities: \ms
\begin{cor}  \it If the Riordan matrix $ \Big(d_{n, k}\Big) = \big(g(x), F(x)\big) $ of finite order has eigenvector  with generating function $v(x)$ then its $n^{\rm th}$ row satisfies
\[\sum_{k = 0}^\infty d_{n, k}\cdot v_k  = \sum_{k = 0}^n d_{n, k}\cdot v_k = g_0\omega^k \cdot v_n.\qquad \Box\]
\end{cor}
\begin{ex}(We thank  L. Woodson for making this formula more elegant.) \[\sum_{k = 0}^{n - 1}  (-1)^k{n \choose k} {2k \choose k}\cdot \frac{1}{4^k} \, = \begin{cases} \hspace{.6in} 0 & \ 
\text{if $n$ is even}\\
&\\
  {2n \choose n}\cdot \frac{2}{4^n} &\ \text{if $n$ is odd}.
  \end{cases} \] 
 \end{ex} 
\ms
{\bf Proof:}
Let the involution $P^\ast$ be gotten by multiplying the odd columns of Pascal's triangle by $(-1)$:
\[P^\ast = \gFpair = \left(\frac{1}{1 - x}, \, \frac{-x}{1 - x}\right)  = \left((-1)^k{n \choose k} \right)_{n,\,  k \geq 0}.\]
 Consider the eigenvector with generating function $v(x) = \sqrt{g(x)} = (1-x)^{-\frac{1}{2}}$: 
\beqn
\sqrt{g(x)} & = &(1 - x)^{-\frac{1}{2}} = \sum_{k = 0}^\infty {-1/2 \choose k}\big(-x)^k\\
& = & \sum_{k = 0}^\infty \frac{1}{k!}\cdot\left(\frac{-1}{2}\cdot\frac{-3}{2}\cdots\frac{-2k + 1}{2}\right)(-1)^kx^k =  \sum_{k = 0}^\infty \frac{1}{k!}\cdot\left(\frac{1}{2}\cdot\frac{3}{2}\cdots\frac{2k -1}{2}\right)x^k \\
& = &  \sum_{k = 0}^\infty \frac{1}{k!}\frac{1}{2^k}\cdot\left(\frac{1\cdot3\cdots (2k - 1))2^kk!}{2^kk!}\right)x^k  = \sum_{k = 0}^\infty\frac{(2k)!}{k!\,k!} \cdot\frac{1}{4^k}\\
& = & \sum_{k = 0}^\infty {2k \choose k}\frac{1}{4^k}\,x^k.
\eeqn
In this example, the eigenvalue in Theorem 3 is $g_0\omega^k = 1\cdot(-1)^0 = 1$. Thus, from Corollary 5.5, the dot product of the $n^{\rm th}$ row with the eigenvector gives 
\[\sum_{k = 0}^n (-1)^k{n\choose k}\cdot{2k \choose k}\frac{1}{4^k} = {2n \choose n}\,\frac{1}{4^n},\]
from which our result follows.\qquad $\Box$
\section {Connection of This Paper to the Papers of Cheon and Kim}
We describe in this section how the results above evolved from our reading of the noted Cheon-Kim papers [4], [5]. We use the notation above in discussing these papers.  Also, we deal only with formal power series, while these authors allow the series to be either formal power series or analytic functions.
\ms
A key concept in their papers is that of an {\sl anti-symmetric series} $\Phi(x, z)$ of two variables and more generally a {\sl $k$-cyclic symmetric series} $\varphi(z_1, \ldots , z_k)$.
\ms
\subsection{Riordan involutions and anti-symmetric series}
The key theorem in [4], relating to L. Shapiro's question (Q9.1,  [16]), is the following:
\bs
\begin{thm*} {\rm ([4], Theorem 2.3 (Q9.1))}
\ms
If $D = \gFpair$ is a Riordan involution then 
\[g(x) = \pm {\rm exp}\left[\Phi(x, F(x)\right] \]
for some anti-symmetric function $\Phi(x, z)$. Conversely, if $F^{(2)}(x) = x \and g(x) = \pm {\rm exp}\left[\Phi(x, F(x)\right]$ for any antisymmetic funcion $\Phi(x, z)$ then $D = \gFpair$ is a Riordan involution.
\end{thm*}\ms
{\bf Comment on the proof of this Theorem:}
\ms
No proof is given in [4].  A reference is given to [12], where this result is stated, also without proof. The proof of the ``Conversely ...''  part of the theorem is straightforward.  But the fact that such a $\Phi(x, z)$ exists, given that $\gFpair$ is a Riordan involution, has not appeared. The Main Results above developed from the author's discovery of the following proof:
\ms
{\bf Proof of Necessity in Cheon-Kim, Theorem 2.3}:
\ms
Let $g(x) = g_0\cdot\wg(x)$.  (Necessarily $g_0 = \pm 1$ since $g(x)\cdot g(F(x)) = 1$.)
Define 
\[\Phi(x, z) = \frac{1}{2}\cdot\ln\left(\frac{\wg(x)}{\wg(z)}\right).\]
Clearly this is antisymmetric: $\Phi(x, z) = -\Phi(z, x)$. Moreover
\beqn
g_0\cdot{\rm exp}\left(\Phi(x, F(x)\right) & = & 
g_0\cdot\exp\left(\frac{1}{2}\cdot\ln\left(\frac{\wg(x)}{\wg(F(x))}\right)\right)\\
& = & g_0\cdot\exp\left(\ln\left(\frac{\wg(x)}{\wg(F(x))}\right)^{\frac{1}{2}}\right)\\
& = & g_0\cdot\exp\left(\ln\left(\frac{\wg(x)}{1/\wg(x)}\right)^{\frac{1}{2}}\right)\\
& = & g_0\cdot\exp\left(\ln(\wg(x))\right) \\
& = & g_0\cdot\wg(x)\\
& = & g(x). \quad \Box
\eeqn
\begin{rmk}\mbox{}
\ms
{\rm  Within this proof we see from the second and the last equations that
\[g(x) = g_0\cdot\frac{\wg^{\frac{1}{2}}(x)}{\wg^{\frac{1}{2}}(F(x))}\]
It was this observation -- that we can let $h(x) = \wg^{\frac{1}{2}}(x)$ - which led to the author's proof of necessity in Theorem 1 in the case $n = 2$. This observation and this proof in this special case led to the explorations which resulted in Theorems 1 and 2 above.  As we shall see below, after we proved Theorem 2 above we saw that it gave the key to proving "the Open question" in [5].
}\end{rmk}
\begin{rmk}
{\rm The ability to use the logarithm and square root (or more generally $n^{\rm th}$ root) functions in these discussions is crucial.  As noted in Remark 1.4, the $n^{\rm th}$ root $g^{\frac{1}{n}}(x)$ is well-defined and the usual laws apply when $g(x) = 1 + G(x)$.
\ms
Similarly, $\ln\big(g(x)\big)$ is defined,  and its laws  justified, in [11] when $g(x) = 1 + G(x)$. We have
\[\ln\big(1 + G(x)\big) \equaldef G(x) - \frac{1}{2}\big(G(x)\big)^2 + \cdots \frac{1}{n}\big(G(x)\big)^n + \cdots \]
\ms
In our definition of $\Phi(x, z)$, application of the logarithm is well-defined only because we may write $\wg(x) = 1 + K(x)$ and we have
\[\frac{\wg(x)}{\wg(z)} = \frac{1 + K(x)}{1+K(z)} = \big(1 + K(x)\big)\cdot\Big(1 - K(z) + \big(K(z)\Big)^2 + \cdots (-1)^n\big(K(z)\Big)^n + \cdots \ \Big)\ =\ 1 + H(x, z).\]
}\end{rmk}
\ms
\begin{rmk}{\rm  {\bf Warning:}
\ms
$\ln(g(x)) = \ln\big(1 + G(x))$ is {\bf not} an element of $\F_1[[x]]$ in the  cases where $g_1 = 0$.  Also, in general as in the proof above, one may not assume that $\Phi(x, F(x))\in \F_1[[x]]$.  So these may {\bf not}, in general,  be taken as the second coordinate of a Riordan element.
\ms
For example, there exist elements $\gFpair$ of every order $n \geq 2$ in the Bell subgroup of ${\cal R}$ for which $\ln(\wg(x)) \notin \F_1[[x]]$.   For we have the following theorem (See [3], Theorem 2.5.4 or [6]):
\begin{thm*}
\ms
If $2\leq n \in \n$\ and $\omega$ is a primitive $n'{\rm th}$\ root of unity in the field $\F$
\ms
{then} for every infinite sequence  $\set{a_k}_{k \neq nj + 1} \in \F$
there exists a
unique sequence \ $\set{a_{nj+1}}_{j=1}^{\infty}$
 such that the formal power series
\ss
\centerline{ $F(x) = \omega x + \sum_{k=2}^{\infty}\ a_k\,x^k$}
\ms
 has order $n$ in $\F_1[[x]]$. \ 
 Moreover, $a_{nj + 1} = 0$ if $a_{2k} = 0$ for all $k \leq j$.
\end{thm*}
\bs
Thus for example, there is an involution $F(x) = - x + a_{10}x^{10 }+ a_{11}x^{11} + \cdots$ with $a_{10}\neq 0$.  Then 
\[ \big(\gFpair\big) = \Big(\frac{F(x)}{x}, F(x)\Big)\]
is an involution in ${\cal R}$, as is $\big(\wg(x), F(x)\big)$ where $\wg(x) = - g(x) = 1 - a_{10}x^9 - \cdots$. It follows that 
$\ln(\wg(x))$ is defined in $\F_+[[x]]$ but it is not an element of $\F_1[[x]]$.  
\ms
This Remark reduces the generality of the validity of Theorem 2.5 and of Corollary 2.6 of [4].  However, the  idea in Corollary 2.6 of [4] of using the logarithm as second coordinate in a conjugating element will be very fruitful in empowering us to give a new proof of C. Marshall's theorem in Section 6.
}\end{rmk}
\subsection {An Answer to an Open Question on $k$-Cyclic Symmetric Series}
\ms
Generalizing the proof above, concerning antisymmetric series and involutions,  we now  answer the ``Open question." in [5] for Riordan elements of higher order.
\begin{defn} 
\ms
A {\bf k-cyclic symmetric series} is a formal series $\varphi(x_1, x_2, \ldots , x_k)$ in $k$ variables such that 
\[\varphi(x_1, x_2, \ldots , x_k) + \varphi(x_2, x_3, \ldots , x_k , x_1) + \cdots \varphi(x_k, x_1, \ldots , x_{k - 1}) = 0\]
\end{defn}
\bs
\begin{qstn*}  ([5], Section 2)

\ms
 If $\gFpair$ is a Riordan element of order $k$ with $g_0 = 1$,  is $g(x)$ of the form\\
$g(x) = \exp\Big(\varphi\big(x, F(x), \ldots , F^{(k - 1)}(x)\big)\Big)$\
for some some $k$-cyclic series   $\varphi(x_1, x_2, \ldots , x_k)${\rm \large ?}
\end{qstn*}
\bs
{\bf Our Approach:}\\
We know from Theorem 1, setting $\wg(x) = g(x)$ since we are given $g_0 = 1$,  that we may write $g(x) = \frac{h(x)}{h\big(F(x)\big)}$.
On the other hand, we may also write
\beqn
g(x) = \exp\left(\ln\big(g(x)\big)\right) & = &\exp\Big(\ln\left(\frac{h(x)}{h\big(F(x)\big)}\right)\Big)\\
\text{so that:} \qquad g(x) = \exp\big(\phi(X)\big) \implies  \phi(X) & = & \ln\left(\frac{h(x)}{h\big(F(x)\big)}\right) \qquad \qquad \qquad \\
\eeqn
However, within the proof of necessity of Theorem 1 above, we see that
\[\frac{h(x)}{h\big(F(x)\big)} = \left[\frac{g(x)^{k - 1}}{g\big(F(x)\big)\cdot  g\big(F^{(2)}(x)\big)\cdot \cdots g\big(F^{(k - 1)}(x)\big) }\right]^{1/k}\]
This leads to the following
\ms
\begin{thm*}[ \bf Affirmative answer to the Open question] \mbox{}
\ms
If $\gFpair$ is a Riordan element of order $k$ with $g_0 = 1$,  then $g(x)$ is of the form
$g(x) = \exp\Big(\varphi\big(x, F(x), \ldots , F^{(k - 1)}(x)\big)\Big)$\
for some some $k$-cyclic series   $\varphi(x_1, x_2, \ldots , x_k)$.
\end{thm*}
\ms
{\bf Proof:}
\ms
Let \ $\varphi(x_1, \ldots , x_k) = \frac{1}{k}\cdot\ln\left[\frac{g(x_1)^{k - 1}}{g(x_2) \ldots g(x_k)}\right]$.
\ms
Then we have:
\be
\item $\varphi(x_1, \ldots , x_k)$ is $k$-cyclic symmetric, by a straightforward calculation using the basic properties of logarithms.
\item $g(x) = \exp\Big(\varphi\big(g(x), g\big(F(x)\big)\ldots , g\big(F^{(k - 1}(x)\big)\big)\Big)$, by the second line of the calculation done in the proof of necessity of Theorem 1 above.
\ee
Therefore the Open question is answered in the affirmative.\qquad $\Box$
\section {A New Proof of C. Marshall's Theorem}  
\ms
C. Marshall's theorem,  in contrast to Theorem 1, deals with the question of finding $F(x)$ for given $g(x)$, which will make  $\gFpair$ into an involution. This theorem shows that there is a unique such $F(x)$ if $g(x)$ is bi-invertible. We shall show how the circle of ideas above can be used to give an alternate (though not as elegant) proof of her Theorem.
\ms
\begin{defn}\mbox{}
\ms
An element $g(x)\in \F[[x]]$ is {\bf bi-invertible} if $g(x) = g_0 + g_1x + g_2x^2  + \cdots = g_0 + G(x)$ where $g_0\neq 0 \and g_1\neq 0$.
\end{defn}
The term ``bi-invertible" refers to the fact that $g(x)$ is invertible in $\F_0[[x]]$ and $G(x)$ is invertible in $\F_1[[x]]$.
\ms
\begin{thm*}{(C. Marshall ([9], [10]))}
\ms
If $g(x)$ is bi-invertible and $\gFpair$ is an involution in ${\cal R}$ then
\[F(x) = \overline{G}\left(\frac{-g_0\cdot G(x)}{g(x)}\right)\]
\end{thm*}
\bs\bs

{\bf Note:} For an involution, necessarily $g_0 = \pm1$. We assume that $F(x) = -x + \cdots$ ({\it i.e., $\omega = -1$)} since the alternative, $\gFpair = (-1, x)$ is not bi-invertible.
\ms
We shall prove this Theorem using the idea in Corollary 2.6 of [4] of conjugating by $\Big(\wg(x)^{\frac{1}{2}}, \ \ln\big(\wg(x)\big)\Big)$. Thus we replace $ \frac{1}{2}\big(x - F(x)\big)$, which (Remark 1.6) always works, by $\ln\big(\wg(x)\big)$.  As pointed out in Remark 5.3, this may be done precisely when $g(x)$ is bi-invertible; {\it i.e.,} when $\ln\big(\wg(x)\big)\in \F_1[[x]]$ .
\ms
{\bf Proof of the Theorem:}  
\be
\item $\big(g(x), F(x)\big) = \Big(\wg(x)^{1/2},\ \ln(\wg(x)\Big)\,\Big(g_0, - x\Big)\,\Big(\wg(x)^{1/2}, \ \ln(\wg(x)\Big)\inverse$.
\beqn
\text{\bf Proof:}\  \gFpair\left(\wg(x)^{\frac{1}{2}}, \ \ln\big(\wg(x)\big)\right) & = & \left(g(x)\cdot\Big(\wg\big(F(x)\big)^{1/2}\Big),\ \ln\Big(\wg\big(F(x)\Big)\right)\\
& = & \left(g_0\cdot\wg(x)\cdot\frac{1}{\wg(x)^{\frac{1}{2}}},\ \ln\Big(\frac{1}{\wg(x)}\Big)  \right)\\
& = & \left(g_0\cdot\wg(x)^{\frac{1}{2}},\ - \ln\big(\wg(x)\big)\right)\\
& = &  \Big(\wg(x)^{1/2},\ \ln(\wg(x)\Big)\Big(g_0, \ -x\Big) \quad \Box
\eeqn
\item We define $L(x)\in \F_1[[x]]$  by $L(x) = \ln(1 + x) = x - \frac{1}{2}x^2 + \frac{1}{3}x^3 + \cdots$. 
\ms
Using the fact that $g_0 = \pm 1$, define 
\[\widehat{G}(x) = \frac{1}{g_0}\cdot G(x) = g_0\cdot G(x) = (\ell_{g_0}\circ G)(x)), \ \text{where} \ \ell_{g_0}(x)\equaldef g_0\cdot x.\] 
 Notice that $\wg(x) = 1 +  \widehat{G}(x)$  and that $\ln(\wg(x)) = \big(L\circ \widehat{G}\big)(x)= \big(L\circ \ell_{g_0}\circ G\big)(x)$.
\ms
  From 1. we have:

\beqn
F(x) & = &\overline{\ln(\wg)}\Big((-id)\big((\ln(\wg)\big)\Big)(x)\\
&&\\
 &=& \overline{\ln(\wg)}\Big(-\ln\big(\wg(x)\big)\Big)\\
 &&\\
& = & \overline{\ln(\wg)}\Big(\ln\Big(\frac{1}{\wg(x)}\Big)\Big)\\
&&\\
& = & \overline{\big(L\circ \ell_{g_0}\circ G\big)}\Big(L\Big(\frac{1}{\wg(x)} - 1\Big)\Big)\\
&&\\
& = & \big(\overline{G}\circ\overline{\ell_{g_0}}\circ\overline{L}\big)\Big(L\Big(\frac{1}{\wg(x)} - 1\Big)\Big)\\
&&\\
& = & \overline{G}\Big(g_0\cdot\frac{1 - \wg(x)}{\wg(x)} \Big) \\
& = & \overline{G}\Big(g_0\cdot\frac{-\widehat{G}(x)}{\wg(x)} \Big) \\
& = & \overline{G}\left(\frac{-g_0\cdot G(x)}{g(x)}\right)\qquad \Box \\
\eeqn
\ee
\bs
\begin{rmk} The uniqueness of $F(x)$ given $g(x)$, indicated by C. Marshall's Theorem, does not hold if the order of $\gFpair$ is greater than two, for the following reason:
\ms
{\bf If} $\gFpair$ has order $n$ and $gcd\big(j, \text{ord}(F)\big) = 1$.
{\bf Then} $\big(g(x),  F^{(j)}(x))$ also has order $n$.\ 
\ms
{\bf Reason}: $1 = g(x)\cdot g(F(x)) \cdots\,\cdot g\big(F^{(n - 1)}x)$\\ 
\mbox{\hspace{.8in}} $= g(x)\cdot g(F^{(j)}(x))\cdot g(F^{(2j)}(x))\cdots g\big(F^{(n - 1)j}(x)\big),$
\ms
as these are exactly the same product, written in different orders. \quad $\Box$
\end{rmk}
\bs 
REFERENCES

\hrulefill
\be
\item{\bf Paul Barry}, {\sl Riordan Arrays: A Primer}, Logic Press (2016).
\item {\bf Salomon Bochner and William Ted Martin}, {\sl Several Complex
Variables}, Princeton University Press (1948).
\item {\bf Thomas S. Brewer}, {\sl Algebraic Properties of Formal Power Series Composition}, University of Kentucky Dissertation (2014).
\item {\bf Gi-Sang Cheon and Hana Kim}, {\sl Simple proofs of open problems about the structure of involutions in the Riordan group}, Linear Algebra and Its Applications 428(2008), 930 - 940.
\item {\bf Gi-Sang Cheon and Hana Kim}, {\sl The Elements of finite order in the Riordan group over the complex field}, Linear Algebra and Its Applications 439(2013), 4032 - 4046.
\item {\bf Marshall M. Cohen}, {\sl Elements of Finite Order in the Group of Formal Power Series Under Composition}, arxiv: 1804.00059 (Presented in MAA Section Meeting April 8, 2006) 
\item {\bf C. Jean-Louis and A. Nkwanta}, {\sl Some Algebraic Structure of the Riordan Group}, Linear Algebra Appl. 438 (2013), 2018-2035.
\item {\bf A. Luz\'{o}n, D. Merlini, M. A. Maron, R. Sprugnoli}, {\sl Identities induced by Riordan Arrays}, Linear Algebra and Its Applications 436 (2012), 631 - 647.
\item {\bf Candice Marshall}, {\sl Another Method of Constructing Pseudo-Involutions in the \break Riordan Group}, Congressus Numerantium, 229 (2017), 343-351.
\item{\bf Candice A. Marshall},{\sl  Construction of Pseudoinvolutions in the Riordan Group}, Morgan State University Dissertation (2017).
\item  {\bf Ivan Niven},  {\sl Formal Power Series}, Amer. Math. Monthly 76 (1969), 871 - 889.
\item {\bf A.D. Polyanin, A.D. Manzhirov},{\sl Handbook of Mathematics for Engineers and Scientists}, Part II Math. Tables: T12.2.1 Functional Equations with Quadratic Nonlinearity, Chapman and Hall, CRC Press, 2006.
\item {\bf Dev Phulara and Louis W. Shapiro}, {\sl Constructing Pseudo-Involutions in the\\ Riordan Group},  Journal of Integer Sequences, Vol. 20 (2017), Article 17.4.7.
\item {\bf Anthony G. O'Farrell and Ian Short}, {\it Reversibility in Dynamics and Group Theory}, London Mathematical Society Lecture Note Series: 416, Cambridge University Press (2015)
\item {\bf L. Shapiro, S. Getu, W. -J. Woan and L. C. Woodson}, {\sl The Riordan Group}, Discrete Applied Mathematics 34 (1991), 229 - 239.
\item {\bf L. W. Shapiro}, {\sl Some open questions about random walks, involutions, limiting distributions and generating functions}, Adv. Appl. Math. 27 (2001) 585-596.
\item {\bf R. Sprugnoli}, {\sl  Riordan Arrays and Combinatorial sums}, Discrete Math. 132 (1994), 267 - 290.
\item {\bf Richard P. Stanley}, {\sl Enumerative Combinatorics, Volume 1, Section 1.1}, Cambridge Studies in Advanced Mathematics 49, Cambridge University Press 1997.
\item  {\bf Richard P. Stanley}, {\sl Enumerative Combinatorics, Volume 2, Section 5.4}, Cambridge Studies in Advanced Mathematics 62, Cambridge University Press 1999.
\ee
\end{document}